# The Diagonalized Newton Algorithm for Non-negative Matrix Factorization

**Hugo Van hamme**
University of Leuven, dept. ESAT
Kasteelpark Arenberg 10 – bus 2441, 3001 Leuven, Belgium
hugo.vanhamme@esat.kuleuven.be

## Abstract

Non-negative matrix factorization (NMF) has become a popular machine learning approach to many problems in text mining, speech and image processing, bio-informatics and seismic data analysis to name a few. In NMF, a matrix of non-negative data is approximated by the low-rank product of two matrices with non-negative entries. In this paper, the approximation quality is measured by the Kullback-Leibler divergence between the data and its low-rank reconstruction. The existence of the simple multiplicative update (MU) algorithm for computing the matrix factors has contributed to the success of NMF. Despite the availability of algorithms showing faster convergence, MU remains popular due to its simplicity. In this paper, a diagonalized Newton algorithm (DNA) is proposed showing faster convergence while the implementation remains simple and suitable for high-rank problems. The DNA algorithm is applied to various publicly available data sets, showing a substantial speed-up on modern hardware.

## 1 Introduction

Non-negative matrix factorization (NMF) denotes the process of factorizing a $N\times T$ data matrix $\mathbf{V}$ of non-negative real numbers into the product of a $N\times R$ matrix $\mathbf{W}$ and a $R\times T$ matrix $\mathbf{H}$, where both $\mathbf{W}$ and $\mathbf{H}$ contain only non-negative real numbers. Taking a column-wise view of the data, i.e. each of the $T$ columns of $\mathbf{V}$ is a sample of $N$-dimensional vector data, the factorization expresses each sample as a (weighted) *addition* of columns of $\mathbf{W}$, which can hence be interpreted as the $R$ parts that make up the data [1]. Hence, NMF can be used to *learn data representations* from samples. In [2], speaker representations are learnt from spectral data using NMF and subsequently applied to separate their signals. Another example in speech processing is [3] and [4], where phone representations are learnt using a convolutional extention of NMF. In [5], time-frequency representations reminiscent of formant traces are learnt from speech using NMF. In [6], NMF is used to learn acoustic representations for words in a vocabulary acquisition and recognition task. Applied to image processing, local features are learnt from examples with NMF in order to represent human faces in a detection task [7].

In this paper, the metric to measure the closeness of reconstruction $\mathbf{Z} = \mathbf{WH}$ to its target $\mathbf{V}$ is measured by their Kullback-Leibler divergence:

$$d_{KL}(\mathbf{V},\mathbf{Z}) = \sum_{n,t} v_{nt} log\left(\frac{v_{nt}}{z_{nt}}\right) - \sum_{n,t} v_{nt} + \sum_{n,t} z_{nt} \tag{1}$$

Given a data matrix **V**, the matrix factors **W** and **H** are then found by minimizing cost function (1), which yields the maximum likelihood estimate if the data are drawn from a Poisson distribution. The multiplicative updates (MU) algorithm proposed in [1] solves exactly this problem in an iterative manner. Its simplicity and the availability of many implementations make it a popular algorithm to date to solve NMF problems. However, there are some drawbacks to the algorithm. Firstly, it only converges locally and is not guaranteed to yield the global minimum of the cost function. It is hence sensitive to the choice of the initial guesses for **W** and **H**. Secondly, MU is very slow to converge. The *goal of this paper* is to *speed up* the convergence while the local convergence property is retained. The resulting Diagonalized Newton Algorithm (DNA) uses only *simple element-wise operations*, such that its implementation requires only a few tens of lines of code, while memory requirements and computational efforts for a single iteration are about the double of an MU update.

The faster convergence rate is obtained by applying Newton's method to minimize $d_{KL}(\mathbf{V},\mathbf{WH})$ over **W** and **H** in alternation. Newton updates have been explored for the Frobenius norm to measure the distance between **V** and **Z** in e.g. [8]-[13]. Specifically, in [11] a diagonal Newton method is applied Frobenius norms. For the Kullback-Leibler divergence, fewer studies are available. Since each optimization problem is multivariate, Newton updates typically imply solving sets of linear equations in each iteration. In [16], the Hessian is reduced by refraining from second order updates for the parameters close to zero. In [17], Newton updates are applied per coordinate, but in a cyclic order, which is troublesome for GPU implementations. In the proposed method, matrix inversion is avoided by diagonalizing the Hessian matrix. The resulting updates resemble the ones derived in [18] to the extent that they involve second order derivatives. Important differences are that [18] involves the non-negative $k$-residuals hence requiring flooring to zero. Of course, the diagonal approximation may affect the convergence rate adversely. Also, Newton algorithms only show (quadratic) convergence when the estimate is sufficiently close to the local minimum and therefore need damping, e.g. Levenberg-Marquardt as in [14], or step size control as in [15] and [16]. In DNA, these convergence issues are addressed by computing both the MU and Newton solutions and selecting the one leading to the greatest reduction in $d_{KL}(\mathbf{V},\mathbf{Z})$. Hence, since the cost is non-decreasing under MU, it will also be under DNA updates. This robust safety net can be constructed fairly efficiently because the quantities required to compute the MU have already been computed in the Newton update. The net result is that DNA iterations are only about two to three times as slow as MU iterations, both on a CPU and on a GPU. The experimental analysis shows that the increased convergence rate generally dominates over the increased cost per iteration such that overall balance is positive and can lead to speedups of up to a factor 6.

## 2 NMF formulation

To induce sparsity on the matrix factors, the KL-divergence is often regularized, i.e. one seeks to minimize:

$$d_{KL}(\mathbf{V},\mathbf{W}\,\mathbf{H}) + \rho \sum_{n,r} w_{nr} + \lambda \sum_{r,t} h_{rt} \tag{2}$$

subject to non-negativity constraints on all entries of **W** and **H**. Here, $\rho$ and $\lambda$ are non-negative regularization parameters.

Minimizing the regularized KL-divergence (2) can be achieved by alternating updates of **W** and **H** for which the cost is non-increasing. The updates for this form of block coordinate descent are:

$$\mathbf{H} \leftarrow \underset{\mathbf{H}' \geq \mathbf{0}}{\arg\min} \left( d_{KL}(\mathbf{V},\mathbf{W}\,\mathbf{H}') + \lambda \sum_{r,t} h'_{rt} \right) \tag{3}$$

$$\mathbf{W} \leftarrow \underset{\mathbf{W}' \geq \mathbf{0}}{\arg\min} \left( d_{KL}(\mathbf{V}, \mathbf{W}'\mathbf{H}) + \rho \sum_{n,r} w'_{nr} \right) \tag{4}$$

Because of the symmetry property $d_{KL}(\mathbf{V},\mathbf{WH}) = d_{KL}(\mathbf{V}^t,\mathbf{H}^t\mathbf{W}^t)$, where superscript-$t$ denotes matrix transpose, it suffices to consider only the update on $\mathbf{H}$. Furthermore, because of the summation over all columns in (1), minimization (3) splits up into $T$ independent optimization problems. Let $\mathbf{v}$ denote any column of $\mathbf{V}$ and let $\mathbf{h}$ denote the corresponding column of $\mathbf{H}$, then the following is the core minimization problem to be considered:

$$\mathbf{h} \leftarrow \underset{\mathbf{h}' \geq \mathbf{0}}{\arg\min} \left( d_{KL}(\mathbf{V}, \mathbf{W}\,\mathbf{h}') + \lambda \mathbf{1}^t \mathbf{h}' \right) \tag{5}$$

where $\mathbf{1}$ denotes a vector of ones of appropriate length. The solution of (5) should satisfy the KKT conditions, i.e. for all $r$ with $h_r > 0$

$$\frac{\partial (d_{KL}(\mathbf{v}, \mathbf{W}\,\mathbf{h}) + \lambda \mathbf{1}^t \mathbf{h})}{\partial h_r} = -\sum_n v_n \frac{w_{nr}}{(\mathbf{W}\,\mathbf{h})_n} + \sum_n w_{nr} + \lambda = 0$$

where $h_r$ denotes the $r$-th component of $\mathbf{h}$. If $h_r = 0$, the partial derivative is positive. Hence the product of $h_r$ and the partial derivative is always zero for a solution of (5), i.e. for $r = 1 \ldots R$:

$$\sum_n v_n \frac{w_{nr} h_r}{(\mathbf{W}\,\mathbf{h})_n} - h_r \left( \sum_n w_{nr} + \lambda \right) = 0 \tag{6}$$

Since $\mathbf{W}$-columns with all-zeros do not contribute to $\mathbf{Z}$, it can be assumed that column sums of $\mathbf{W}$ are non-zero, so the above can be recast as:

$$h_r \frac{(\mathbf{W}^t \mathbf{q})_r}{(\mathbf{W}^t \mathbf{1})_r + \lambda} - h_r = 0$$

where $q_n = v_n \,/\, (\mathbf{Wh})_n$. To facilitate the derivations below, the following notations are introduced:

$$a_r = \frac{(\mathbf{W}^t \mathbf{q})_r}{(\mathbf{W}^t \mathbf{1})_r + \lambda} - 1 \tag{7}$$

which are functions of $\mathbf{h}$ via $\mathbf{q}$. The KKT conditions are hence recast as [20]

$$a_r\, h_r = 0 \quad \text{for } r = 1 \ldots R \tag{8}$$

Finally, summing (6) over $r$ yields

$$\sum_r ((\mathbf{W}^t \mathbf{1})_r + \lambda) h_r = \sum_n v_n \tag{9}$$

which is satisfied for any guess $\mathbf{h}$ by renormalizing:

$$h_r \leftarrow h_r \frac{\mathbf{v}^t \mathbf{1}}{\mathbf{h}^t (\mathbf{W}^t \mathbf{1} + \lambda)} \tag{10}$$

## 2.1 Multiplicative updates

For the more generic class of Bregman divergences, it was shown in a.o. [20] that multiplicative updates (MU) are non-decreasing at each update of $\mathbf{W}$ and $\mathbf{H}$. For KL-divergence, MU are identical to a fixed point update of (6), i.e.

$$h_r \leftarrow h_r (1 + a_r) = h_r \frac{(\mathbf{W}^t \mathbf{q})_r}{(\mathbf{W}^t \mathbf{1})_r + \lambda} \tag{11}$$

Update (11) has two fixed points: $h_r = 0$ and $a_r = 0$. In the former case, the KKT conditions imply that $a_r$ is negative.

## 2.2 Newton updates

To find the stationary points of (2), $R$ equations (8) need to be solved for $\mathbf{h}$. In general, let $\mathbf{g}(\mathbf{h})$ be an $R$-dimensional vector function of an $R$-dimensional variable $\mathbf{h}$. Newton's update then states:

$$\mathbf{h} \leftarrow \mathbf{h} - (\nabla\mathbf{g})^{-1}\mathbf{g}(\mathbf{h}) \quad \text{with} \quad (\nabla\mathbf{g})_{rl} = \frac{\partial g_r(\mathbf{h})}{\partial h_l} \tag{12}$$

Applied to equations (8):

$$(\nabla\mathbf{g})_{rl} = a_l\delta_{rl} - \frac{h_r}{(\mathbf{W}^t\mathbf{1})_r + \lambda}\sum_i \frac{v_i w_{ir} w_{il}}{(\mathbf{W}\mathbf{h})_i^2} \tag{13}$$

where $\delta_{rl}$ is Kronecker's delta. To avoid the matrix inversion in update (12), the last term in equation (13) is diagonalized, which is equivalent to solving the $r$-th equation in (8) for $h_r$ with all other components fixed. With

$$b_r = \frac{1}{(\mathbf{W}^t\mathbf{1})_r + \lambda}\sum_n v_n \frac{w_{nr}^2}{(\mathbf{W}\,\mathbf{h})_n^2} \tag{14}$$

which is always positive, an element-wise Newton update for $\mathbf{h}$ is obtained:

$$h_r \leftarrow h_r \frac{h_r b_r}{h_r b_r - a_r} \tag{15}$$

Notice that this update does not automatically satisfy (9), so updates should be followed by a renormalization (10). One needs to pay attention to the fact that Newton updates will attract towards both local minima and local maxima. Like for the EM-update, $h_r = 0$ and $a_r = 0$ are the only fixed points of update (15), which are now shown to be *locally* stable. In case the optimizer is at $h_r = 0$, $a_r$ is negative by the KKT conditions, and update (15) will indeed decrease $h_r$. In a sufficiently small neighborhood of a point where the gradient vanishes, i.e. $a_r = 0$, update (15) will increase (decrease) $h_r$ if and only if (11) increases (decreases) its estimate. Since if (11) never increases the cost, update (15) attracts to a minimum.

However, this only guarantees local convergence for per-element updates and Newton methods are known to suffer from potentially small convergence regions. This also applies to update (15), which can indeed result in limit cycles in some cases. In the next subsections, two measures are taken to respectively increase the convergence region and to make the update non-increasing.

## 2.3 Step size limitation

When $a_r$ is positive, update (15) may not be well-behaved in the sense that its denominator can become negative or zero. To respect nonnegativity and to avoid the singularity, it is bounded below by a function with the same local behavior around zero:

$$\frac{h_r b_r}{h_r b_r - a_r} = \frac{1}{1 - \frac{a_r}{h_r b_r}} \geq 1 + \frac{a_r}{h_r b_r} \tag{16}$$

Hence, if $a_r \geq 0$, the following update is used:

$$h_r \leftarrow h_r\left(1 + \frac{a_r}{h_r b_r}\right) = h_r + \frac{a_r}{b_r} \tag{17}$$

Finally, step sizes are further limited by flooring resp. ceiling the multiplicative gain applied to $h_r$ in update (15) and (17) (see Algorithm 1, steps 11 and 24 for details).

## 2.4 Non-increase of the cost

Despite the measures taken in section 2.3, the divergence can still increase under the Newton update. A very safe option is to compute the EM update additionally and compare the cost function value for both updates. If the EM update is be better, the Newton update is rejected and the EM update is taken instead. This will guarantee non-increase of the cost function. The computational cost of this operation is dominated by evaluating the KL-divergence, not in computing the update itself.

## 3 The Diagonalized Newton Algorithm for KLD-NMF

In Algorithm 1, the arguments given above are joined to form the Diagonalized Newton Algorithm (DNA) for NMF with Kullback-Leibler divergence cost. Matlab™ code is available from www.esat.kuleuven.be/psi/spraak/downloads both for the case when **V** is sparse or dense.

Algorithm 1: pseudocode for the DNA KLD-NMF algorithm. $\oslash$ and $\odot$ are element-wise division and multiplication respectively and $[x]_\varepsilon = \max(x,\varepsilon)$. Steps not labelsed with "MU" is the additional code required for DNA.

Input: data **V**, initial guess for **W** and **H**, regularization weights $\rho$ and $\lambda$.
MU - Step 1: divide the $r$-th column of **W** by $\sum_n w_{nr} + \lambda$. Multiply the $r$-th row of **H** by the same number.
MU - Step 2: **Z** = **WH**
MU - Step 3: $\mathbf{Q} = \mathbf{V} \oslash \mathbf{Z}$
**Repeat until convergence**
    MU - Step 4: precompute $\mathbf{W} \odot \mathbf{W}$
    MU - Step 5: $\mathbf{A} = \mathbf{W}^t\mathbf{Q} - 1$
    MU - Step 6: $\mathbf{H}_{MU} = \mathbf{H} + \mathbf{A} \odot \mathbf{H}$
    MU - Step 7: $\mathbf{Z}_{MU} = \mathbf{W}\mathbf{H}_{MU}$
    MU - Step 8: $\mathbf{Q}_{\boldsymbol{MU}} = \mathbf{V} \oslash \mathbf{Z}_{MU}$
    MU - Step 9: $\mathbf{d}_{MU} = \mathbf{1}^t(\mathbf{V} \odot \log(\mathbf{Q}_{MU}))$
    Step 10: $\bar{\mathbf{Q}} = \mathbf{V} \oslash (\boldsymbol{Z} \odot \boldsymbol{Z})$; $\mathbf{B} = (\mathbf{W} \odot \mathbf{W})^t\bar{\mathbf{Q}}$
    Step 11: $\mathbf{H}_{DNA} = \mathbf{H} \odot [\mathbf{B} \oslash (\mathbf{B} - \mathbf{A})]_{\boldsymbol{\varepsilon}}$ for the entries for which **A** < 0
        $\mathbf{H}_{DNA} = \mathbf{H} + \mathbf{min}(\mathbf{A} \oslash (\mathbf{H} \odot \mathbf{B}), \alpha\mathbf{H})$ for the entries for which **A** ≥ 0
        multiply $t$-th column of $\mathbf{H}_{DNA}$ with the $t$-th entry of $(\mathbf{1}^t\mathbf{V}) \oslash (\mathbf{1}^t\mathbf{H}_{DNA})$
    Step 12: $\mathbf{Z}_{DNA} = \mathbf{W}\mathbf{H}_{DNA}$
    Step 13: $\mathbf{Q}_{\boldsymbol{DNA}} = \mathbf{V} \oslash \mathbf{Z}_{\boldsymbol{DNA}}$
    Step 14: $\mathbf{d}_{DNA} = \mathbf{1}^t(\mathbf{V} \odot \log(\mathbf{Q}_{DNA}))$
    Step 15: copy **H**, **Z** and **Q** from:
        $\mathbf{H}_{DNA}$, $\mathbf{Z}_{DNA}$ and $\mathbf{Q}_{DNA}$ for the columns for which $\mathbf{d}_{DNA} < \mathbf{d}_{MU}$
        $\mathbf{H}_{EM}$, $\mathbf{Z}_{EM}$ and $\mathbf{Q}_{EM}$ for the columns for which $\mathbf{d}_{DNA} \geq \mathbf{d}_{MU}$
    MU - Step 16: divide (multiply) the $r$-th row (column) of **H** (**W**) by $\sum_t h_{nt} + \rho$.
    Step 17: precompute $\mathbf{H} \odot \mathbf{H}$
    MU - Step 18: $\mathbf{A} = \mathbf{Q}\mathbf{H}^{\boldsymbol{t}} - 1$
    MU - Step 19: $\mathbf{W}_{MU} = \mathbf{W} + \mathbf{A} \odot \mathbf{W}$
    MU - Step 20: $\mathbf{Z}_{MU} = \mathbf{W}_{MU}\mathbf{H}$
    MU - Step 21: $\mathbf{Q}_{\boldsymbol{MU}} = \mathbf{V} \oslash \mathbf{Z}_{MU}$
    MU - Step 22: $\mathbf{d}_{MU} = (\mathbf{V} \odot \log(\mathbf{Q}_{MU}))\mathbf{1}$
    Step 23: $\bar{\mathbf{Q}} = \mathbf{V} \oslash (\boldsymbol{Z} \odot \boldsymbol{Z})$; $\mathbf{B} = \bar{\mathbf{Q}}(\mathbf{H} \odot \mathbf{H})^{\boldsymbol{t}}$
    Step 24: $\mathbf{W}_{DNA} = \mathbf{W} \odot [\mathbf{B} \oslash (\mathbf{B} - \mathbf{A})]_{\boldsymbol{\varepsilon}}$ for the entries for which **A** < 0
        $\mathbf{W}_{DNA} = \mathbf{W} + \mathbf{min}(\mathbf{A} \oslash (\mathbf{H} \odot \mathbf{B}), \alpha\mathbf{W})$ for the entries for which **A** ≥ 0
        multiply the $n$-th row of $\mathbf{W}_{DNA}$ with the $n$-th entry of $(\mathbf{V1}) \oslash (\mathbf{W}_{DNA}\mathbf{1})$
    Step 25: $\mathbf{Z}_{DNA} = \mathbf{W}_{DNA}\mathbf{H}$
    Step 26: $\mathbf{Q}_{\boldsymbol{DNA}} = \mathbf{V} \oslash \mathbf{Z}_{\boldsymbol{DNA}}$
    Step 27: $\mathbf{d}_{DNA} = (\mathbf{V} \odot \log(\mathbf{Q}_{DNA}))\mathbf{1}$
    Step 28: copy **W**, **Z** and **Q** from:
        $\mathbf{W}_{DNA}$, $\mathbf{Z}_{DNA}$ and $\mathbf{Q}_{DNA}$ for the rows for which $\mathbf{d}_{DNA} < \mathbf{d}_{MU}$
        $\mathbf{W}_{EM}$, $\mathbf{Z}_{EM}$ and $\mathbf{Q}_{EM}$ for the rows for which $\mathbf{d}_{DNA} \geq \mathbf{d}_{MU}$
    MU - Step 29: divide(multiply) the $r$-th column (row) of **W** (**H**) by $\sum_n w_{nr} + \lambda$.

Notice that step 9, 14 22 and 27 require some care for the zeros in **V**, which should not contribute to the cost. In terms of complexity, the most expensive steps are the computation of **A**, **B**, $\mathbf{Z}_{MU}$ and $\mathbf{Z}_{DNA}$, which require O(*NRT*) operations. All other steps require O(*NR*), O(*RT*) or O(*NR*) operations. Hence, it is expected that a DNA iteration is about twice as slow as MU iteration. On modern hardware, parallelization may however distort this picture and hence experimental verification is requied.

# 4 Experiments

DNA and MU are run on several publicly available[1] data sets. In all cases, **W** is initialized with a random matrix with uniform distribution, normalized column-wise. Then **H** is initialized as $\mathbf{W}^t\mathbf{V}$ and one MU iteration is performed. The same initial values are used for both algorithms. Sparsity is not included in this study, so $\rho = \lambda = 0$. The algorithm parameters are set to $\varepsilon = 0.01$ and $\alpha$=4. CPU timing measurements are obtained on a quad-core AMD™ Opteron 8356 processor running the MATLAB™ code available at www.esat.kuleuven.be/psi/spraak/downloads which uses the built-in parallelization capability. Timing measurements on the graphical processing unit (GPU) are obtained on a TESLA C2070 running MATLAB and Accelereyes Jacket v2.3.

## 4.1 Dense data matrices

The first dataset considered is a set of 400 frontal face greyscale 64×64 images of 40 people showing 10 different expressions. The resulting 4096×165 dense matrix is decomposed with factors of a common dimension *R* of 10, 20, 40 and 80. Figure 1 shows the KL divergence as a function of iteration number and CPU time as measured on the CPU. The superiority of DNA is obvious: for instance, at *R* = 40, DNA reaches the same divergence after 33 iterations as MU obtains after 500 iterations. This implies a speed-up of a factor 15 in terms of iterations or 6.3 in terms of CPU time.

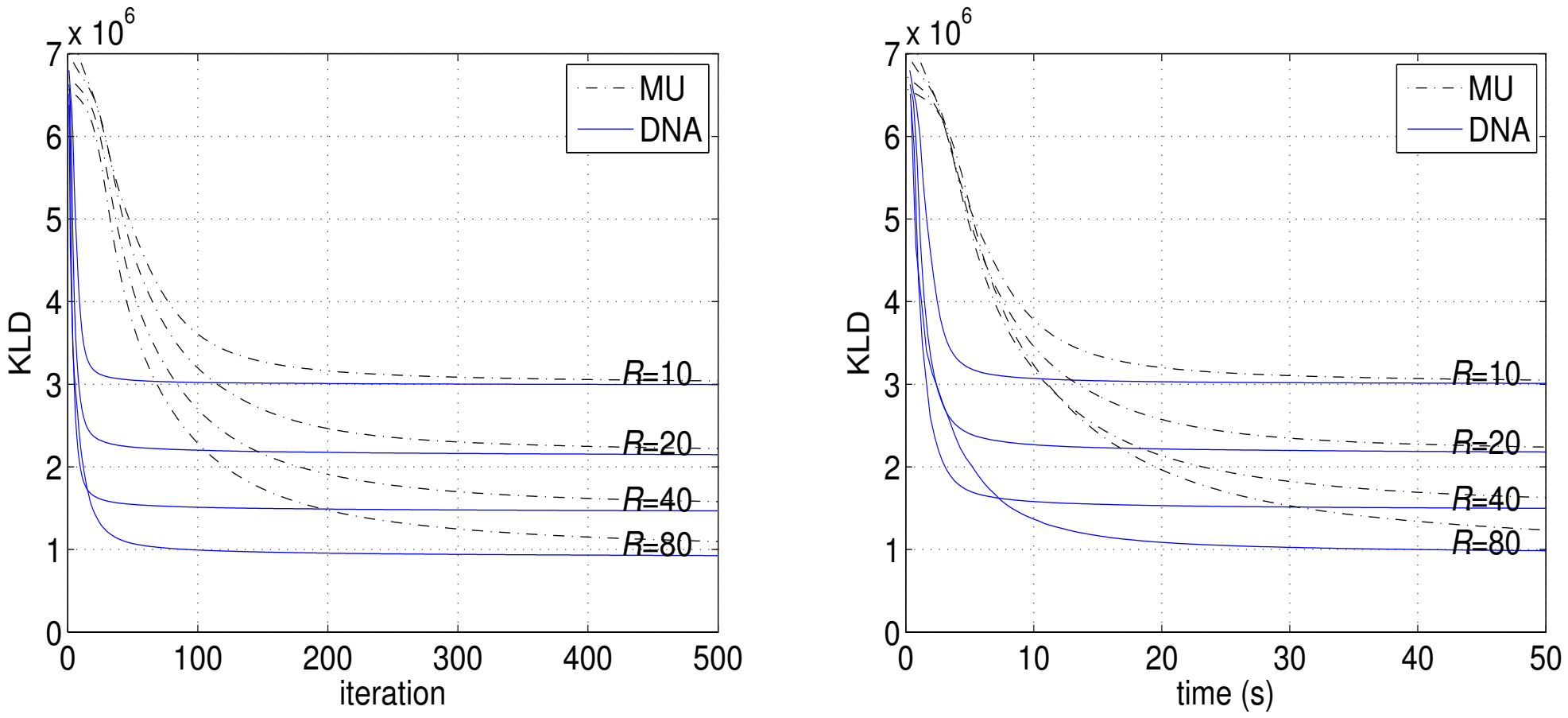

Figure 1: convergence of DNA and MU on the ORL image dataset as a function of the number of iterations (left) and CPU time (right) for different ranks *R*.

The second test case is the CMU PIE dataset which consists of 11554 greyscale images of 32×32 pixels showing human faces under different illumination conditions and poses. The data are shaped to a dense 1024×11554 matrix and a decomposition of rank *R* = 10, 20, 40 and 80 are attempted with the MU and DNA algorithms. As observed in Figure 2, the proposed DNA still outperforms MU, but by a smaller margin.

[1] www.cad.zju.edu.cn/home/dengcai/Data/data.html

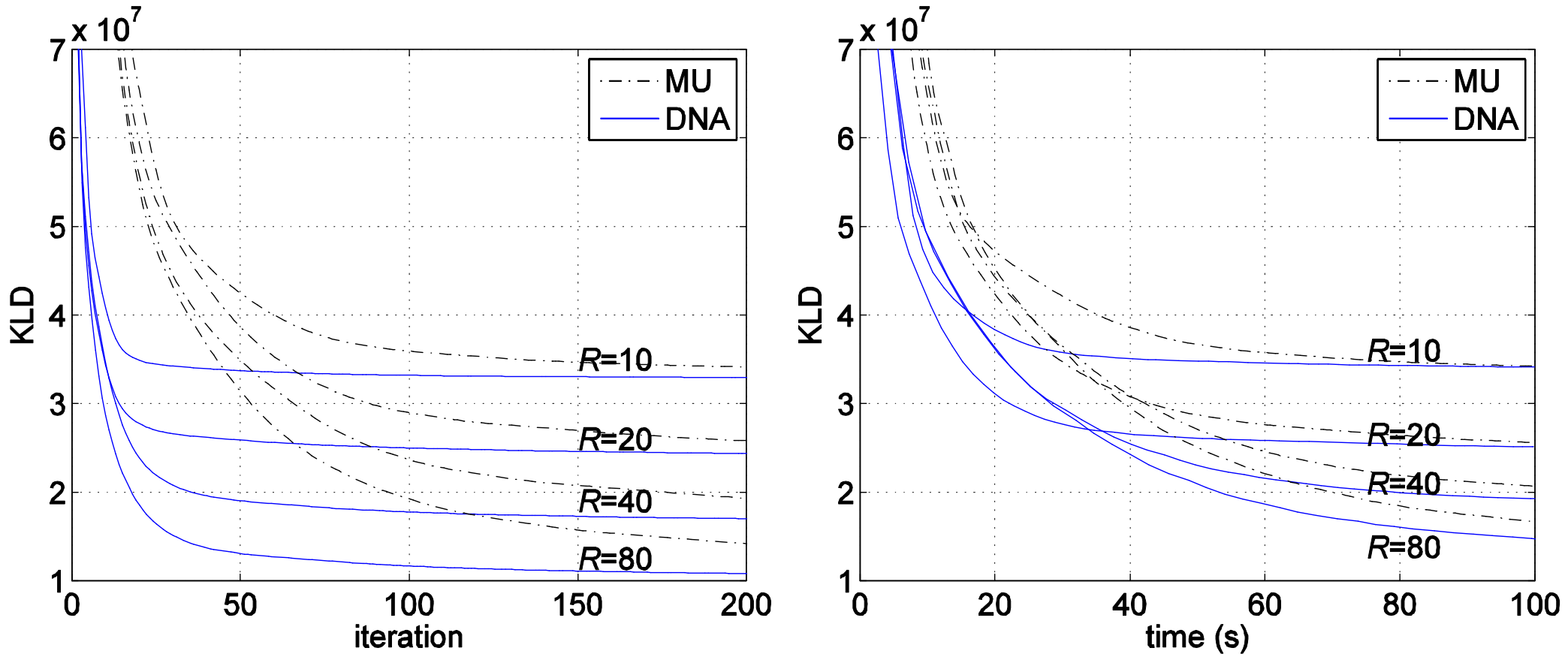

Figure 2: convergence of DNA and MU on the CMU PIE image dataset as a function of the number of iterations (left) and CPU time (right).

An overview of the time required for a single iteration on both data sets is given in Table 1. For MU, the first row lists the time if the KL divergence is not computed as this is not required if the number of iterations is fixed in advance instead of stopping the algorithm based on a decrease in KLD. The table shows that the computational cost of MU can be reduced by about a third by not computing KLD. Compared to MU with cost calculation, DNA requires typically about 2.5 to 3 times more time per iteration on the CPU. On the GPU, the ratio is rather 2 to 2.5.

### 4.2 Sparse data matrices

The third matrix considered originates from the NIST Topic Detection and Tracking Corpus (TDT2). For 10212 documents (columns of **V**), the frequency of 36771 terms (rows of **V**) was counted leading to a sparse $36771 \times 10212$ matrix with only 0.35% non-zero entries. The fourth matrix originates from the Newsgroup corpus results in a 61188×18774 sparse frequency matrix with 0.2% non-zeros. Both for MU and DNA a MATLAB implementation using the sparse matrix class was made. In this case, an iteration of DNA is twice as slow a MU iteration. Again, the convergence of both algorithms is shown in Figure 3. In this case, DNA is only marginally faster than MU in terms of CPU time.

Table 1: time per iteration in milliseconds as measured on the CPU and GPU implementations for different ranks ($R$) and dense matrices (ORL/PIE).

| dataset | ORL | | | | | | | | PIE | | | | | | | |
|---|---|---|---|---|---|---|---|---|---|---|---|---|---|---|---|---|
| R | 10 | | 20 | | 40 | | 80 | | 10 | | 20 | | 40 | | 80 | |
| processor | CPU | GPU | CPU | GPU | CPU | GPU | CPU | GPU | CPU | GPU | CPU | GPU | CPU | GPU | CPU | GPU |
| MU without cost | 78 | 3.7 | 85 | 4.1 | 96 | 5.0 | 115 | 6.8 | 310 | 18 | 310 | 20 | 330 | 27 | 400 | 38 |
| MU with cost | 114 | 6.4 | 118 | 7.0 | 130 | 7.9 | 161 | 9.4 | 480 | 35 | 490 | 37 | 510 | 44 | 580 | 55 |
| DNA | 269 | 15.5 | 280 | 15.9 | 319 | 17.9 | 425 | 23.6 | 1180 | 71 | 1430 | 76 | 1720 | 95 | 1960 | 125 |

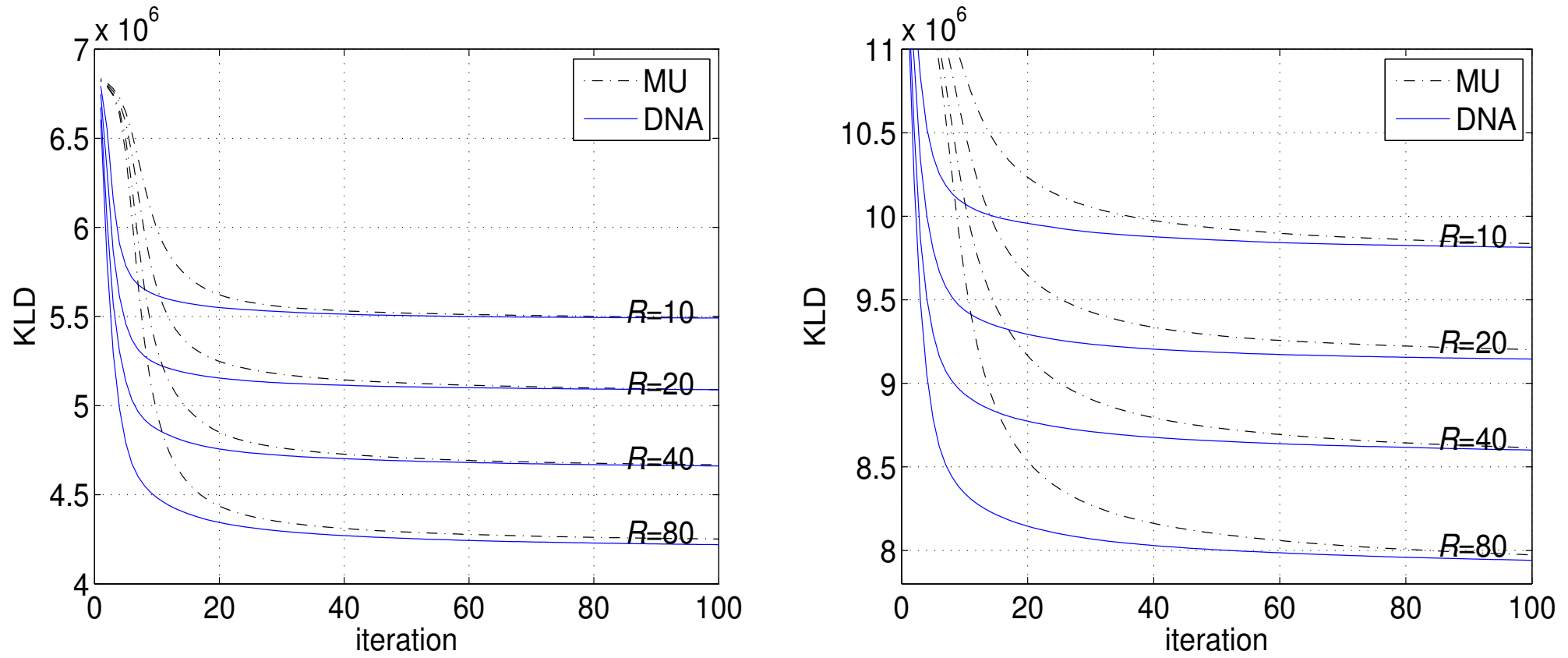

Figure 3: convergence of DNA and MU on the sparse TDT2 (left) and Newsgroup (right) data.

## 5 Conclusions

The DNA algorithm is based on Newton's method for solving the stationarity conditions of the constrained optimization problem implied by NMF. This paper only addresses the Kullback-Leibler divergence as a cost function. To avoid matrix inversion, a diagonal approximation is made, resulting in element-wise updates. Experimental verification on publicly available matrices with a CPU and GPU MATLAB implementation for dense data matrices and a CPU MATLAB implementation for sparse data matrices show that, depending on the case and matrix sizes, DNA iterations are 2 to 3 times slower than MU iterations. In most cases, the diagonal approximation is good enough such that faster convergence is observed and a net gain results.

Since Newton updates can in general not ensure monotonic decrease of the cost function, the step size was controlled with a brute force strategy of falling back to MU in case the cost is increased. More refined step damping methods could speed up DNA by avoiding evaluations of the cost function, which is next on the research agenda.

### Acknowledgments

This work is supported by IWT-SBO project 100049 (ALADIN) and by KU Leuven research grant OT/09/028(VASI).

## References

[1] D. Lee, and H. Seung, "Algorithms for non-negative matrix factorization," *Advances in Neural Information Processing Systems*, vol. 13, pp. 556–562, 2001.

[2] B. Raj, R. Singh and P. Smaragdis, "Recognizing Speech from Simultaneous Speakers", *in proceedings of Eurospeech*, pp. 3317-3320, Lisbon, Portugal, September 2005

[3] P. Smaragdis, "Convolutive Speech Bases and Their Application to Supervised Speaker Separation," *IEEE Transactions on Audio, Speech and Language Processing*, vol. 15, pp. 1-12, January 2007

[4] P. D. O'Grady and B. A. Pearlmutter, "Discovering Speech Phones Using Convolutive Non-negative Matrix Factorisation with a Sparseness Constraint.," *Neurocomputing*, vol. 72, no. 1-3, pp. 88-101, December 2008, ISSN 0925-2312.

[5] M. Van Segbroeck and H. Van hamme, "Unsupervised learning of time-frequency patches as a noise-robust representation of speech," *Speech Communication*, volume 51, no. 11, pp. 1124-1138, November 2009.

[6] H. Van hamme, "HAC-models: a Novel Approach to Continuous Speech Recognition," In Proc. International Conference on Spoken Language Processing, pp. 2554-2557, Brisbane, Australia, September 2008.

[7] X. Chen, L. Gu, S. Z. Li and H.-J. Zhang, "Learning representative local features for face detection," *in proceedings of the IEEE Computer Society Conference on Computer Vision and Pattern Recognition*, pp. 1126-1131, Kauai, HI, USA, December 2001.

[8] D. Kim, S. Sra and I. S. Dhillon, "Fast Projection-Based Methods for the Least Squares Nonnegative Matrix Approximation Problem," *Statistical Analy Data Mining*, vol. 1, 2008

[9] C.-J. Lin, "Projected gradient methods for non-negative matrix factorization," *Neural Computation*, vol. 19, pp. 2756-2779, 2007

[10] R. Zdunek, A. H. Phan and A. Cichocki, "Damped Newton Iterations for Nonnegative Matrix Factorization," *Australian Journal of Intelligent Information Processing Systems*, 12(1), pp. 16-22, 2010

[11] Y. Zheng, and Q. Zhang, "Damped Newton based Iterative Non-negative Matrix Factorization for Intelligent Wood Defects Detection," *Journal of software*, vol. 5, no. 8, pp. 899-906, August 2010.

[12] P. Gong, and C. Zhang, "Efficient Nonnegative Matrix Factorization via projected Newton method", *Pattern Recognition*, vol. 45, no. 9, pp. 3557-3565, September 2012.

[13] S. Bellavia, M. Macconi, and B. Morini, "An interior point Newton-like method for nonnegative least-squares problems with degenerate solution," *Numerical Linear Algebra with Applications*, vol. 13, no. 10, pp. 825-846, December 2006.

[14] R. Zdunek and A. Cichocki, "Non-Negative Matrix Factorization with Quasi-Newton Optimization," *Lecture Notes in Computer Science, Artificial Intelligence and Soft Computing* 4029, pp. 870-879, 2006

[15] R. Zdunek and A. Cichocki, "Nonnegative Matrix Factorization with Constrained Second-Order Optimization", *Signal Processing*, vol. 87, pp. 1904-1916, 2007

[16] G. Landi and E. Loli Piccolomini, "A projected Newton-CG method for nonnegative astronomical image deblurring," *Numerical Algorithms*, no. 48, pp. 279–300, 2008

[17] C.-J. Hsieh and I. S. Dhillon, "Fast Coordinate Descent Methods with Variable Selection for Non-negative Matrix Factorization," *in proceedings of the 17th ACM SIGKDD International Conference on Knowledge Discovery & Data Mining* (KDD), San Diego, CA, USA, August 2011

[18] L. Li, G. Lebanon and H. Park, "Fast Bregman Divergence NMF using Taylor Expansion and Coordinate Descent," *in proceedings of the 18th ACM SIGKDD Conference on Knowledge Discovery and Data Mining,* 2012

[19] A. Cichocki, S. Cruces, and S.-I. Amari, "Generalized Alpha-Beta Divergences and Their Application to Robust Nonnegative Matrix Factorization," *Entropy*, vol. 13, pp. 134-170, 2011; doi:10.3390/e13010134

[20] I. S. Dhillon and S. Sra, "Generalized Nonnegative Matrix Approximations with Bregman Divergences," *Neural Information Proc. Systems*, pp. 283-290, 2005